\documentclass[11pt]{amsart}
\usepackage{amssymb,amscd,amsfonts}
\usepackage[english]{babel}
\usepackage[latin1]{inputenc}
\usepackage[OT1]{fontenc}
\usepackage{times}
\usepackage{amsmath}
\usepackage{amsthm}
\usepackage{amsfonts}
\usepackage{amssymb}
\usepackage[all]{xy}
\newtheorem{theorem}{Theorem}[section]

\newtheorem{proposition}[theorem]{Proposition}

\theoremstyle{remark}
\newtheorem{remark}[theorem]{Remark}

\newtheorem{example}[theorem]{Example}

\newcommand{\Lie}[1]{\mathfrak{#1}}
\newcommand{\lLie}[1]{\mathfrak{#1}^{\vee}}
\newcommand{\dLie}[1]{\mathfrak{#1}^{\ast}}
\newcommand{\Weyl}{\mathcal{W}}
\newcommand{\lWeyl}{\mathcal{W}^{\vee}}
\newcommand{\Sg}{S(\mathfrak{g}^*)}
\newcommand{\Sh}{S(\mathfrak{h}^*)}
\newcommand{\lSh}{S((\mathfrak{h}^\vee)^*)}
\newcommand{\Sgg}{S(\mathfrak{g}^*)^\mathfrak{g}}
\newcommand{\lSgg}{S((\mathfrak{g}^\vee)^*)^{\mathfrak{g}^\vee}}
\newcommand{\ShW}{S(\mathfrak{h}^*)^{\mathcal{W}}}
\newcommand{\lShW}{S((\mathfrak{h}^\vee)^*)^{\mathcal{W}^\vee}}
\newcommand{\C}{\mathbf{C}}
\newcommand{\Ch}{\text{Chev}}
\newcommand{\np}{\mathfrak{n}_+}

\newcommand{\rank}{\text{rank}}

\newcommand{\Ker}{\text{Ker}}
\newcommand{\ad}[1]{\text{ad}_{#1}}

\newcommand{\Int}{\text{Int}}
\newcommand{\sltwo}{\textit{sl}_2}

\begin{document}

\title{Principal basis in Cartan subalgebra}

\author{Rudolf Philippe Rohr}
\address{University of Geneva, Department of Mathematics,
2-4 rue du Lièvre, c.p. 64, 1211 Geneva 4, Switzerland}
\email{Rudolf.Rohr@math.unige.ch}

\date{\today}

\begin{abstract} Let $\Lie{g}$ be a simple complex Lie algebra and $\Lie{h}$ a Cartan subalgebra. In this article we explain how to obtain the principal basis of $\Lie{h}$ starting form a set of generators $\{p_1, \cdots ,p_r\}$,$r=\rank(\Lie{g})$, of the invariants polynomials $\Sgg$. For each invariant polynomial $p$, we define a $G$-equivariant map $Dp$ form $\Lie{g}$ to $\Lie{g}$. We show that the Gram-Schmidt orthogonalization of the elements  $\{Dp_1(\rho^\vee), \cdots Dp_r(\rho^\vee) \}$ gives the principal basis of $\Lie{h}$. Similarly the orthogonalization of the elements $\{Dp_1(\rho), \cdots Dp_r(\rho) \}$ produces the principal basis of the Cartan subalgebra of $\lLie{g}$, the Langlands dual of $\Lie{g}$.

\end{abstract}

\maketitle
\tableofcontents

%---------------------------------------------------------------------------------------------

\section{Introduction}
\label{sec_intro}
Let $\Lie{g}$ be a simple Lie algebras over $\C$, let $\Lie{h}$ be a choice of the Cartan subalgebra, and let $\{\alpha_1, \cdots , \alpha_r\}$ be a choice for the set of simple roots.. We consider the principal $\sltwo$-triple, $\Lie{s}_0=<e_0,f_0,h_0>_\C$, given by (Lemma 5.2 of \cite{Ko1})\\

\begin{tabular}{ll}
	the semi simple element : & $h_0= 2 \rho^{\vee}=$the sum of positive coroots \\
	the positive nilpotent elements : & $e_0=e_1 + \cdots + e_r=$the sum of positive nilpotent \\
	& element corresponding to simple roots \\
	the negative nilpotent element : & $f_0$ such that the well know commutation relation \\
	& are satisfied, i.e. $[h_0,e_0]=2e_0, [h_0,f_0]=-2f_0, [e_0,f_0]=h_0.$
\end{tabular}\\

The restriction of the adjoint action of $\Lie{g}$ to $\Lie{s}_0$ gives a representation of $\Lie{s}_0$ on $\Lie{g}$. Let 

\begin{equation}
\label{equ_deci}
	\Lie{g} = V_{1} \oplus \cdots \oplus V_{r}
\end{equation}
be the decomposition of this representation into irreducible $\Lie{s}_0$-modules. There is exactly $r=\rank(\Lie{g})$ irreducible submodules and their dimensions are $\dim(V_i)=2k_i+1$, where the $\{k_i\}_{i=1}^{n}$ are the exponents of $\Lie{g}$ (\cite{Ko1}, Section 5). Without loss of generality, we can suppose $1=k_1 \leq k_2 \leq \cdots \leq k_r$. Each of the modules $V_i$ has a one dimensional intersection with $\Lie{h}$, i.e. for each $i$ there exists $h_i$ such that $<h_i>_\C = V_{i} \cap \Lie{h}$. By definition (\cite{Ba}, Section 7.2) these elements form the principal basis of $\Lie{h}$, 

$$\text{principal basis}=\{h_1, \cdots h_r\}.$$

The principal basis was introduce by Kostant to solve a question about Clifford algebra (see remark \ref{rem_conj}).\\

If all exponents are distinct, the principal basis is uniquely defined (up to constant multiples). This is the case for all simples Lie algebras except for $so(2l,\C)$ with $l$ even. In this case there are exactly two exponents with the value $l-1$. Then in the decomposition (\ref{equ_deci}) there are two modules of the same dimension. Hence the definition of the principal basis has to be refined in this case (see §\ref{sec_d2l}).\\

Let $\{\tilde{p}_1, \cdots ,\tilde{p}_r\}$ be a choice of homogeneous generators in the ring of invariant polynomials $\Sgg$. We suppose $\deg(p_1) \leq \cdots \leq \deg(p_r)$. Let $\{e^*_i\}_i$ be an orthonormal basis of $\dLie{g}$ with respect to a choice of an invariant bilinear form $B$. We define 

$$d\tilde{p}_i(x) = \sum_j \frac{\partial{\tilde{p}_i}}{\partial e^*_i}(x) \otimes e^*_i \qquad x \in \Lie{g}$$
as a $G$-equivariant map $d\tilde{p}_i: \Lie{g} \rightarrow \dLie{g}$, and we consider the set 

$$\{d\tilde{p}_1(\rho^\vee), \cdots , d\tilde{p}_r(\rho^\vee) \}$$
of elements of $\dLie{g}$. Using $B$ we can identify $\Lie{g}$ with $\dLie{g}$ and we can consider these elements as the elements in $\Lie{g}$. Our first result is (Theorem \ref{thm_main_1})

\begin{enumerate}
	\item The elements $\{d\tilde{p}_1(\rho^\vee), \cdots , d\tilde{p}_r(\rho^\vee) \}$ form a basis of $\Lie{h}$.
	\item After orthogonalization we obtain the principal basis of $\Lie{h} \subset \Lie{g}$.
\end{enumerate}

Let $\rho$ be the half sum of positive roots. Using $B$ to identify $\Lie{g}$ with $\dLie{g}$ we construct the set 

$$\{d\tilde{p}_1(\rho), \cdots , d\tilde{p}_r(\rho) \}$$
of elements of $\dLie{g}$. Our second result is (Theorem \ref{thm_main_2})

\begin{enumerate}
	\item The elements $\{d\tilde{p}_1(\rho), \cdots , d\tilde{p}_r(\rho) \}$ form a basis of $\lLie{h}$, the Cartan subalgebra of $\lLie{g}$, the Langlands dual of $\Lie{g}$.
	\item After orthogonalization we obtain the principal basis of $\lLie{h} \subset \lLie{g}$.\\
\end{enumerate}

In §\ref{sec_dual} we discuss the relation between roots ans coroots and between the rings $\Sgg$, $\ShW$, $\lSgg$ and $\lShW$.\\

In §\ref{sec_tds} we define the principal three dimensional subalgebra and the principal basis.\\

§\ref{sec_conj} contains the proof of the first two results, Theorem \ref{thm_main_1} and \ref{thm_main_2}.\\

In §\ref{sec_d2l} we study in detail the case of $so(2l,\C)$ ($l$ even). This Lie algebra has two exponents with the same value, and the principal basis defined by the decomposition (\ref{equ_deci}) is not unique. Let $\sigma$ be the automorphism of order two of the Dinkin diagram of type $D_l$ ($l$ even). If we require that the vectors of the principal basis be eigenvectors of $\sigma$, then the principal basis is unique (up to constant multiple) and orthogonal (Theorem \ref{thm_auto}).\\

\textbf{Acknowledgments} I would like to thank Anton Alekseev for very helpful discussions. This work was supported in part by the Swiss National Science Foundation.

%---------------------------------------------------------------------------------------------

\section{Dual root system}
\label{sec_dual}

Throughout this section $\Lie{g}$ denotes a simple Lie algebra over the field of complex numbers, $\Lie{h} \subset \Lie{g}$ a Cartan subalgebra, $\mathcal{R}$ the root system and $B$ a non-degenerate invariant bilinear form on $\Lie{g}$, e.g. the Killing form or the canonical form. \\

The bilinear form $B$ induces an isomorphism of vectors spaces between $\Lie{g}$ and $\dLie{g}$, 

$$B^\natural : \Lie{g} \stackrel{\cong}{\longrightarrow} \dLie{g},$$
and its inverse 

$$B^\flat : \dLie{g} \stackrel{\cong}{\longrightarrow} \Lie{g}.$$
They are given by $B^\natural(x)=B(x,\cdot)$ and $B^\flat(\alpha)=x_\alpha$, where $x_\alpha$ is the unique element of $\Lie{g}$ such that $B(x_\alpha,\cdot)=\alpha$. As the restriction of $B$ to the Cartan subalgebra $\Lie{h}$ is non degenerate, the restriction of $B^\natural$ and $B^\flat$ to the Cartan subalgebra are isomorphisms. Moreover $B$ induces a non degenerate bilinear form on $\Lie{h}^\ast$ which will by denoted again by $B$ and it is given for all $\alpha,\beta \in \dLie{h}$ by 

$$B(\alpha,\beta)=B(B^\flat(\alpha),B^\flat(\beta)).$$\\

To each root $\alpha \in \dLie{h}$ we associate a coroot $\alpha^\vee = \gamma(\alpha) \in \Lie{h}$,

\begin{equation}
\label{equ_root}
	\Lie{h}^* \ni \alpha \stackrel{\gamma}{\longrightarrow} \alpha^\vee = \frac{2}{B(\alpha,\alpha)}B^\flat(\alpha) \in \Lie{h}.
\end{equation}
The dual roots play the role of the roots in the dual root system $\mathcal{R}^\vee$ of $\Lie{g}$. The map $\gamma$ defines an isomorphism in the sense that $\gamma : \mathcal{R} \stackrel{\cong}{\longrightarrow} \mathcal{R}^\vee$. Note that $\mathcal{R}^\vee \subset \Lie{h}$. The map $\gamma$ is canonical, i.e. it is independent of the bilinear form $B$, but it does not extend to a vector space isomorphism between $\dLie{h}$ and $\Lie{h}$.\\ 

Now we will recall the definition of the Langlands dual $\lLie{g}$ of $\Lie{g}$.

\textit{
\begin{center}
The Langlands dual $\lLie{g}$ is the Lie algebra whose root system is $\mathcal{R}^\vee$.
\end{center}}
This definition is equivalent to :

\textit{
\begin{center}
The Cartan matrix of $\lLie{g}$ is the transpose of the Cartan matrix of $\Lie{g}$.
\end{center}}
or

\textit{
\begin{center}
The Dynkin diagram of $\lLie{g}$ is that of $\Lie{g}$ with arrow reversed.
\end{center}}
The Langlands dual of the Langlands dual of $\Lie{g}$ is $\Lie{g}$. For the root systems of simple Lie algebras we have the $X_n^\vee = X_n$ for $X=A,D,E,F,G$ and $B_n^\vee = C_n$.\\

As $\mathcal{R}$ spans $\Lie{h}^*$, dual to the Cartan subalgebra of $\Lie{g}$, the dual root system $\mathcal{R}^\vee$ spans $(\lLie{h})^*$, dual to the Cartan subalgebra of $\lLie{g}$. We have the following identifications 

$$\lLie{h} = \dLie{h} \stackrel{\cong \text{ by }B}{\longleftrightarrow} \Lie{h} = (\lLie{h})^*.$$
Note that the left and right identifications are canonical, but the middle one depends on $B$.\\

Let $\Weyl$ and $\lWeyl$ be the Weyl groups of $\Lie{g}$ and $\lLie{g}$ respectively. Let $\sigma_{\alpha} \in \Weyl$ and $\sigma_{\alpha^\vee} \in \lWeyl$ be reflections by root $\alpha$ and coroot $\alpha^\vee$ respectively. We have that $B^\flat$ intertwine them, i.e.

$$ B^\flat \circ \sigma_{\alpha} = \sigma_{\alpha^\vee} \circ B^\flat.$$
The isomorphism $B^\flat$ extends to a graded algebra isomorphism

$$B^\flat :\Sh \stackrel{\cong}{\longrightarrow} \lSh,$$
and because it intertwines the action of the Weyl groups, its restriction to the Weyl group invariants is a graded algebra isomorphism

$$B^\flat :\ShW  \stackrel{\cong}{\longrightarrow} \lShW.$$\\

Using the famous result \cite{Ch} of Chevalley we have that $\ShW$ and $\lShW$ are generated by $r=\rank(\Lie{g})$ homogeneous linearly independent polynomials. Let $\{p_1 , \cdots p_r\}$ be a choice of generators of $\ShW$. Then $\{p^\flat_1=B^\flat(p_1),\cdots,p^\flat_r=B^\flat(p_r)\} $ are generators of $\lShW$, i.e.

\begin{equation}
\label{equ_gen}
	\ShW = \C[p_1, \cdots , p_r] \stackrel{\cong}{\rightarrow} \lShW=\C[p^\flat_1,\cdots,p^\flat_r],
\end{equation}
Moreover the $p_i$'s (and also the $p^\flat_i$'s) are homogeneous of degree $k_i+1$, where the integers $\{k_i\}_{i=1}^{r}$ are the exponents of $\Lie{g}$. Note that $\Lie{g}$ and $\lLie{g}$ have the same exponents. For simple Lie algebra we can choice these generators so that $1=k_1 \leq k_2 \leq \cdots \leq k_r$ and the $p_1$ is given by the restriction of $B$ to $\dLie{h}$. 

\begin{remark}
Using the Tables of \cite{Bo}, we notice that for a simple Lie algebra its exponents are all different, except for simple Lie algebras with roots system of type $D_l$ with $l$ even and greater or equal to $4$.
\end{remark}

Let $\Ch : \Sg \rightarrow \Sh$ be the Chevalley projection. Its restriction to the invariants induces a graded algebra isomorphism (see \cite{Di}, Chapter 7, Section 3), i.e $\Ch : \Sgg \stackrel{\cong}{\rightarrow} \ShW$. This isomorphism implies that $\Sgg$ and $\lSgg$ are generated by $r$ linearly independent polynomials. We obtain the following commutative diagram of graded algebra isomorphisms :

\begin{equation}
\label{equ_diag_iso}
	\xymatrix{
		 \Sgg =\C[\tilde{p}_1,\cdots,\tilde{p}_r] \ar[d]_{\Ch}  \ar[r]^{\cong}  & \lSgg =\C[\hat{p}_1,\cdots,\hat{p}_r] \ar[d]^{\Ch} \\
		 \ShW = \C[p_1, \cdots , p_r] \ar[r]^{B^\flat} & \lShW = \C[p^\flat_1, \cdots , p^\flat_r]}
\end{equation}
Moreover we can choose the generators so that : 

\begin{enumerate}
	\item $\Ch(\tilde{p}_i) = p_i$ and $\Ch(\hat{p}_i) = p^\flat_i,$
	\item the upper horizontal isomorphism sends the generators $\tilde{p}_i$ to the generators $\hat{p}_i$.
\end{enumerate}
Note that $\tilde{p}_1$ (resp. $\hat{p}_1$) are defined up to a constant by an invariant bilinear form of $\Lie{g}$ (resp. $\dLie{g}$).\\

%---------------------------------------------------------------------------------------------

\section{Decomposition under action of the principal three dimensional subalgebra, and the principal basis of the Cartan subalgebras}
\label{sec_tds}

Let $\Lie{s}=<e,f,h>_\C$ be a $\sltwo$-triple, \footnote{The commutation relations are $[h,e]=2e$, $[h,f]=-2f$ and $[e,f]=h$.} of a simple complex Lie algebra $\Lie{g}$. The restriction of the adjoint representation of $\Lie{g}$ to $\Lie{s}$ give a representation of $\Lie{s}$ on $\Lie{g}$, i.e.  

$$ \Lie{s} \ni x \rightarrow \ad{x} \in End(\Lie{g}) .$$ 
is a Lie algebra homomorphism. Note that $\sltwo$-triple exist, indeed every nilpotent elements of $\Lie{g}$ can by a Morosov's Theorem embedded in a TDS (see Section 3 of \cite{Ko1}). We consider the decomposition of $\Lie{g}$ into a direct sum of the irreducible $\Lie{s}$-modules

\begin{equation}
\label{equ_dec}
	\Lie{g} = V_1 \oplus V_2 \oplus \dots \oplus V_n.
\end{equation}
Without loss of generality we suppose that $\dim(V_1) \leq \dim(V_2) \leq \cdots \leq \dim(V_n)$. In general the number $n$ of irreducible $\Lie{s}$-modules is greater that the rank of $\Lie{g}$ (see \cite{Ko1}, Section 5), and if its equal to $r$ then $\Lie{s}$ will be call principal $\sltwo$-triple (Theorem 5.2 in \cite{Ko1}).\\

We define a particular principal $\sltwo$-triple, $\Lie{s}_0 = <h_0,e_0,f_0>_\C$, with

\begin{equation}
\label{equ_pTDS}
	h_0 = 2 \rho^\vee , \qquad  e_0 = \sum_i e_i \quad \text{and} \quad f_0 = \sum_i c_i f_i,
\end{equation}
where $\rho^\vee$ is the half sum of positive coroots, $e_i$ (rep. $f_i$) are the root vectors corresponding to the simple roots $\alpha_i$ (resp. $-\alpha_i$) such that $B(e_i,f_i)=1$ and the $c_i$ are given by the relation $\sum_i c_i B^\flat(\alpha_i) = 2\rho^\vee$. Moreover each principal $\sltwo$-triple is conjugate to (\ref{equ_pTDS}), and if we require the semisimple element to be in $\Lie{h}$, then the principal $\sltwo$-triple is given by (\ref{equ_pTDS}). It is of course possible to multiply each $e_i$ by a scalar and then divide $f_i$ by the same scalar.\\

Using $\Lie{s}_0$ in the decomposition (\ref{equ_dec}) we have that $n=\rank(\Lie{g})$ and each $V_i$ is an irreducible $\Lie{s}_0$-module. Moreover $\dim(V_i)=2k_i+1$, where the $\{k_i\}_{i=1}^{r}$ are the exponents of $\Lie{g}$. This facts implies that for all $i$ the intersection of $V_i$ and $\Lie{h}$ is one dimensional, i.e. 

\begin{center}
\textit{for all $i$, there exist $h_i \in \Lie{h}$ such that $<h_i>_\C = V_i \cap \Lie{h}$.} \\
\end{center}

These elements form a basis of $\Lie{h}$. This basis is called the principal basis of $\Lie{h}$ (\cite{Ba}, Section 7), 

\begin{equation}
\label{equ_pBasis}
	\text{principal basis of }\Lie{h} = \{h_1, \cdots ,h_r \}.
\end{equation}

For $k_i \neq k_j$ we have that $h_i$ is orthogonal to $h_j$ relatively to $B$. Indeed, we have that $ad_{e_0} ad_{f_0} (h_i) = \frac{k_i(k_i+1)}{2} h_i$, then using the invariance of $B$ we conclude that for $k_i \neq k_j$, $h_i$ is orthogonal to $h_j$. Hence if all exponents are distinct, up to a constant, the principal basis is unique. In fact a similar argument shows that the decomposition (\ref{equ_dec}) under $\Lie{s}_0$ is an orthogonal decomposition. The only cases when the are two equal exponent occur for the Lie algebras with the root system of type $D_{l}$ with $l$ even, i.e. for $so(2l,\C)$. We study this case in §\ref{sec_d2l}. 

\begin{remark}
	As vector space $V_1 = \Lie{s}_0$, then we have (up to a constant) $h_1=\rho^\vee$.
\end{remark}

\begin{remark}
	The elements $\{h_i\}_{i=1}^r$  of the principal basis are characterize by 
	\begin{enumerate}
		\item $h_i \in \Ker(ad^{k_i+1}(e_0))$, 
		\item $h_i \notin \Ker(ad^{k_i}(e_0))$,
		\item	for $i \neq j$, $h_i \bot h_j$.
	\end{enumerate}
\end{remark}

\begin{remark}
	If an element $h \in \Lie{h}$ satisfies $ad^{k_i+1}(e_0)(h)=0$, then $h \in (V_{1} \oplus \cdots \oplus V_{i}) \cap \Lie{h}$.
\end{remark}

We will end this section by giving two examples of principal basis.

\begin{example}
\label{ex_g2}
Let $\Lie{g}$ be the Lie algebra of type $G_2$. Its Cartan subalgebra $\Lie{h}$ is $2$-dimensional. Let $\alpha_1$ and $\alpha_2$ be a choice of the simples roots. We choose $B$ to be the canonical form, i.e. on simple root it is given by
$$B(\alpha_1,\alpha_1)=2 \qquad B(\alpha_1,\alpha_2)= -1 \qquad B(\alpha_2,\alpha_2)=\frac{2}{3}.$$
The coroots are then given by 
$$\alpha_1^{\vee}= B^\flat(\alpha_1) \qquad \alpha_2^{\vee} = 3 B^\flat(\alpha_2) .$$
The first vector in the principal basis is the sum of the positive coroots, and the second is one orthogonal to the first. They are (up to constant) 
$$h_1 = 5 \alpha_1^{\vee} + 3 \alpha_2^{\vee}  \qquad h_2 = -3 \alpha_1^{\vee} + \alpha_2^{\vee}.$$
\end{example}

\begin{example}[Principal basis for $\text{sl}(3,\C)$]
\label{ex_sl3}
Let $e_{ij}$ be the canonical basis of $\text{Mat}(3,\C)$. The standard choice of basis for the Cartan subalgebra of $\text{sl}(3,\C)$ is $\tilde{h}_1=e_{11}-e_{22}$ and $\tilde{h}_2=e_{22}-e_{33}$. The positive (resp. negative) nilpotent elements are given by $e_{12}$, $e_{23}$ and $e_{13}$ (resp. $e_{21}$, $e_{32}$ and $e_{31}$). In terms of this basis the principal $\sltwo$-triple $\Lie{s}_0$ is given by $h_0= 2\tilde{h}_1+2\tilde{h}_2$, $e_0=e_{12}+e_{23}$ and $f_0=2e_{21}+2e_{32}$. In the decomposition (\ref{equ_dec}) the two irreducible $\Lie{s}_0$-modules are given by
$$V_1=<e_0,h_0,f_0>_\C,$$  
$$V_2 = <e_{13}, e_{23}-e_{12}, \tilde{h}_{1}-\tilde{h}_{2}, e_{32}-e_{21},e_{31}>_\C.$$
Then the principal basis is given by (up to constants)
$$\{h_1=\tilde{h}_1+\tilde{h}_2=e_{11}-e_{33},h_2=\tilde{h}_1-\tilde{h}_2=e_{11}-2e_{22}+e_{33} \}.$$
By a direct calculation we verify that it is an orthogonal basis.\\
\end{example}

%---------------------------------------------------------------------------------------------

\section{How to obtain the principal basis of the Cartan subalgebra using the generators of the invariants polynomials}
\label{sec_conj}

Let $\Lie{g}$ be a simple Lie algebra over the complex field, and $\tilde{p} \in \Sgg$ an invariant polynomial. We consider $d\tilde{p}$ the differential of this polynomial. It is defined by 

$$d\tilde{p} = \sum_i \frac{\partial \tilde{p}}{\partial e^*_i} \otimes e^*_i \in \Sg \otimes \dLie{g},$$
where $\{e^*_i\}$ is an orthonormal basis of $\dLie{g}$. We view this differential as an map form $\Lie{g}$ to $\dLie{g}$, i.e.

$$\Lie{g} \ni x \longrightarrow d\tilde{p}(x) = \sum_i \frac{\partial \tilde{p}}{\partial e^*_i}(x) \cdot e^*_i \in \dLie{g}.$$

Let $G= \Int(\Lie{g})$ be the group of inner automorphisms of $\Lie{g}$. It is generated by the elements of the form $e^{\ad{x}}$ with $x \in \Lie{g}$ nilpotent. The invariant bilinear form $B$ is invariant by $G$, i.e. for all $g \in G$ and $x,y \in \Lie{g}$ we have $B(g \cdot x ,y)=B(x,g^{-1} \cdot y)$. Moreover, $G$ intertwines the isomorphisms $B^\flat$ and $B^\natural$. 

\begin{proposition}
\label{prop_equi}
The map $d\tilde{p}$ is $G$-equivariant, i.e. for all $g \in G$

$$d\tilde{p} \circ g = g \circ d\tilde{p}.$$
\end{proposition}
\begin{proof}
Let $x \in \Lie{g}$ and $g \in G$, then $d\tilde{p}(g \cdot x)= \sum_i \frac{\partial \tilde{p}}{e_i}( g \cdot x)e_i$. Let $y_i = g^{-1}x_i$ be a new orthonormal basis of $\Lie{g}$. Then using the fact that $p$ is $G$-invariant we obtain

$$\sum_i \frac{\partial \tilde{p}}{e_i}( g \cdot x)e_i = \sum_i \frac{\partial \tilde{p}}{g \cdot y_i}( g \cdot x)g \cdot y_i = \sum_i \frac{\partial (g \cdot \tilde{p})}{g \cdot y_i}( g \cdot x)g \cdot y_i = \sum_i \frac{\partial \tilde{p}}{y_i}(x)g \cdot y_i.$$
\end{proof}

Let $D\tilde{p}$ be the composition of $d\tilde{p}$ with $B^\flat$, i.e. 
$$D\tilde{p} = B^\flat \circ d\tilde{p} : \Lie{g} \longrightarrow \Lie{g}.$$
Obviously $D\tilde{p}$ is $G$-equivariant.\\

We can consider the same construction with an invariant polynomial $p \in \ShW$. Then we obtain the map

$$Dp :\Lie{h} \longrightarrow \Lie{h}.$$

The following proposition relates these two constructions.

\begin{proposition}
\label{prop_link}
Let $p \in \ShW$ and $\tilde{p} \in \Sgg$ be invariant polynomials such that $\Ch(\tilde{p})=p$. Then for all $h \in \Lie{h}$ we have

$$Dp(h)=D\tilde{p}(h).$$
\end{proposition}
\begin{proof}
	By Lemma 7.3.6 of \cite{Di}, there exist $a \in \Sgg \np^*$ such that $\tilde{p} = p + a$. But $Da(h)=0$ because $\dLie{h}$ is orthogonal to $\np^*$.
\end{proof}

Let $\{\tilde{p}_1, \cdots ,\tilde{p}_r\}$ be a choice of generators of $\Sgg$, with $\deg(p_1) \leq \deg(p_2) \leq \cdots \leq \deg(p_r)$. Let $\rho^\vee \in \Lie{h}$ be the half sum of positive coroots. By the preceding proposition the elements 
\begin{equation}
\label{equ_base_avant}
	\{D\tilde{p}_1(\rho^{\vee}),\cdots D\tilde{p}_r(\rho^\vee) \} 
\end{equation}
are in $\Lie{h}$. Note that up to a constant we can chose $\tilde{p}_1$ such that $D\tilde{p}_1(\rho^{\vee})=\rho^{\vee}$. By Theorem 3 of \cite{Va} with $A=\rho^\vee$, this family is linearly independent. Hence it is a basis of $\Lie{h}$. Using the Gram-Schmidt orthogonalization process on this basis we obtain the orthogonal basis :

\begin{equation}
\label{equ_basis_after}
\xymatrix{
\{D\tilde{p}_1(\rho^{\vee}),\cdots D\tilde{p}_r(\rho^\vee) \} \ar[d]^{\text{Gram Schmidt}}\\
\\
{
\begin{split}
	\{ & h_1=D\tilde{p}_1(\rho^{\vee}) , \\
	& h_2=D\tilde{p}_2(\rho^{\vee}) - \lambda_{21} D\tilde{p}_1(\rho^{\vee}), \\ 
	& \cdots , \\
	& h_r=D\tilde{p}_r(\rho^\vee) - \lambda_{r,r-1} D\tilde{p}_1(\rho^{\vee}) - \cdots \}.
\end{split}
}
}
\end{equation}

The first result of this section is the following Theorem. 

\begin{theorem}
\label{thm_main_1}
		The orthogonal basis (\ref{equ_basis_after}) of $\Lie{h}$ is (up to constant multiple) the principal basis of the Cartan subalgebra $\Lie{h}$ of $\Lie{g}$.\\
\end{theorem}

Before we prove this Theorem we will show how it work in the Examples \ref{ex_g2} and \ref{ex_sl3}.

\begin{example}[Continuation of Example \ref{ex_g2}]
We will compute the principal basis using Theorem \ref{thm_main_1} and then compare with the previous result. Let $\{ x,y \}$ be the orthonormal basis of the Cartan subalgebra given by 
$$x = \frac{\alpha_1^\vee}{\sqrt{2}} \qquad y = \sqrt{\frac{2}{3}}\left(\alpha_2^\vee + \frac{3}{2} \alpha_1^\vee \right).$$
Let $\{x^*,y^*\}$ be the dual basis. A choice for generators of $\ShW$ is given by
$$p_1 = (x^*)^2 + (y^*)^2 \qquad p_2 = 33(x^*)^6+27(y^*)^6+45(x^*)^4(y^*)^2 +135(x^*)^2(y^*)^4.$$
The differentials of these polynomials computed at $\rho^{\vee}$ give (up to constant multiple)
$$Dp_1{\rho^{\vee}} = \text{cte} \rho^{\vee} \qquad Dp_2{\rho^\vee} = \text{cte} 2425 \alpha_1^\vee + 1383 \alpha_2^\vee.$$
We remark that they are not orthogonal, bur after orthogonalization we obtain the principal basis. 
\end{example}

\begin{remark}
The previous example shows that the orthogonalization is indispensable.
\end{remark}

\begin{example}[Continuation of example \ref{ex_sl3}]
We will compute the principal basis using Theorem \ref{thm_main_1} and then compare with the previous result. A choice of generators of $\ShW$ is given by (\cite{Me}, Section 3),
$$p_1 = (e^*_{11})^2+(e^*_{22})^2+(e^*_{33})^2  \qquad p_2 = (e^*_{11})^3+(e^*_{22})^3+(e^*_{33})^3 ,$$
where $\{e^*_{11},e^*_{22},e^*_{33} \}$ is the dual basis of $\{e_{11},e_{22},e_{33} \}$. As explain in \cite{Me}, these polynomials belong to the dual of the Cartan subalgebra of $gl(3,\C)$. We work with the basis $\{ \tilde{h}^*_1=e^*_{11}-e^*_{22} , \tilde{h}^*_2=e^*_{22}-e^*_{33} \}$ completed by the central element $c^* = e^*_{11} + e^*_{22} + e^*_{33}$. In this basis these polynomials are written as follows: 
$$p_1 =(\tilde{h}^*_1)^2 + (\tilde{h}^*_2)^2 + \tilde{h}^*_1\tilde{h}^*_2 + (c^*)^2 \qquad p_2 = 2(\tilde{h}^*_1)^3  - 2 (\tilde{h}^*_2)^3 + 3 (\tilde{h}^*_1)^2 \tilde{h}^*_2 - 3 \tilde{h}^*_1 (\tilde{h}^*_2)^2 + \text{ terms with }c^*.$$
Remark that $<\tilde{h}^*_1,\tilde{h}^*_2>_\C \bot <c^*>_\C$. Their differentials evaluated at $\rho^{\vee}$ are (up to a constant)
$$Dp_1(\rho^{\vee}) = \tilde{h}_1 + \tilde{h}_2 \qquad Dp_2(\rho^{\vee}) = \tilde{h}_1 - \tilde{h}_2.$$
This is fortunately the same basis as before. Note that in this case the choice of polynomial is unique (up to constants) and we do not need to orthogonalize.
\end{example}

\begin{proof}[Proof of Theorem \ref{thm_main_1}]
This proof work if all exponents are different.  This is the case for all simple Lie algebras except for type $D_{l}$ with $l$ even. We proof this case in §\ref{sec_d2l}.\\

\begin{enumerate}
	\item Using Theorem 3 of \cite{Va} with $A=\rho^\vee$, we conclude that the elements (\ref{equ_base_avant}) are linearly independents.
	\item Let $\Lie{s}_0$ be the principal $\sltwo$-triple given by (\ref{equ_pTDS}). Then $D\tilde{p}_i(\rho^\vee) \in \Ker(\ad{e_0}^{k_i+1})$. Indeed using the $G$-equivariance of $D\tilde{p}_i)$ (Proposition \ref{prop_equi}) for $e^{t\ad{e_0}}$ ($t$ a parameter) we have 
	$$(e^{t\ad{e_0}} D\tilde{p}_i)(\rho^\vee) = D\tilde{p}_i(e^{t\ad{e_0}}\rho^\vee) = (D\tilde{p}_i \circ(1+t\ad{e_0}))(\rho^\vee),$$
	but and $D\tilde{p}_i$ is a polynomial of degree $k_i$, and the term in $t^{k_i+1}$ is vanishes.\\
	
	\item The step (b) implies that for all $i$ 
	$$D\tilde{p}_i(\rho^\vee) \in \left( V_{1} \oplus \cdots \oplus V_{i} \right) \cap \Lie{h} .$$
	
	\item Then if we orthogonalize the family 
	$$\{D\tilde{p}_1(\rho^\vee),\cdots,D\tilde{p}_r(\rho^\vee) \}$$
	as (\ref{equ_basis_after}) we obtain (up to constant multiple) the principal basis of $\Lie{h}$.
\end{enumerate}
\end{proof}

This Theorem provides a method to compute the principal basis of the Cartan subalgebra of $\Lie{h}$ using a set of generators of the invariants polynomials $\Sgg$. Now we will give the second result of this section, which explain how to obtain the principal basis of the Cartan subalgebra of $\lLie{g}$ the Langlands dual of $\Lie{g}$ using a set of generators of $\Sgg$.\\

As explained before $d\tilde{p}$ is an $G$-equivariant map from $\Lie{g}$ to $\dLie{g}$. We consider the following map, which is also $G$-equivariant,

$$\hat{D}\tilde{p} = d\tilde{p} \circ B^\natural : \dLie{g} \longrightarrow \dLie{g}.$$ 

We have the analogue of Proposition \ref{prop_link}

\begin{proposition}
\label{prop_link2}
Let $p \in \ShW$ and $\tilde{p} \in \Sgg$ be invariant polynomials such that $\Ch(\tilde{p})=p$. Then for all $\lambda \in \dLie{h}$ we have

$$\hat{D}p(\lambda)=\hat{D}\tilde{p}(\lambda).$$
\end{proposition}
The proof is the same as for Proposition \ref{prop_link}.\\

The second result of this section is the following theorem.

\begin{theorem}
\label{thm_main_2}
	Let $\rho$ be the half sum of the positive roots of $\Lie{g}$. Let $\{\hat{D}\tilde{p_1}(\rho), \cdots , \hat{D}\tilde{p_r}(\rho) \}$ be a family of element of $\dLie{h}$, the Cartan subalgebra of $\lLie{g}$. Then
	\begin{enumerate}
		\item This family is linearly independent.
		\item	Using the Gram-Schmidt process as in (\ref{equ_basis_after}) we obtain (up to constant) the principal basis of $\lLie{h} \subset \lLie{g}$.
		\end{enumerate}
\end{theorem}
\begin{proof}
	\begin{enumerate}
		\item Let $p \in \ShW$ be an invariant polynomial and define $p^\flat = B^\flat(p) \in \lShW$. As for $p$, we consider $dp^\flat$ to be a map from $\dLie{h}$ to $\Lie{h}$. We have the following relation between $dp$ and $dp^\flat$,
		
		$$dp \circ B^\flat = B^\natural \circ dp^\flat.$$
		\item We take $\tilde{p}\in \Sgg$ such that $\Ch(\tilde{p})=p$. Using the Proposition \ref{prop_link2} and (a) we have for all $\lambda \in \dLie{h}$,
		
		$$\hat{D}\tilde{p}(\lambda)=\hat{D}p(\lambda)=(B^\natural \circ dp^\flat)(\lambda).$$
		\item We take $\hat{p} \in \lSgg$ such that $\Ch(\hat{p})=p^\flat$. We use the equalities in (b) to obtain
		
		$$\hat{D}\tilde{p}(\lambda)=\hat{D}p(\lambda)=(B^\natural \circ dp^\flat)(\lambda)=(B^\natural \circ d\hat{p})(\lambda).$$
		The proof of the last equality is the same as for the Proposition $\ref{prop_link}$ with $\lLie{g}$ at the place of $\Lie{g}$.
		\item
		By (c) point we see that the family  $\{\hat{D}\tilde{p_1}(\rho), \cdots , \hat{D}\tilde{p_r}(\rho) \}$ is the same as the family $\{(B^\natural \circ d\hat{p_1}^\flat)(\rho), \cdots ,(B^\natural \circ d\hat{p_r}^\flat)(\rho) \}$. The result (a) is given by Theorem 3 in \cite{Va} and the result (b) by using Theorem \ref{thm_main_1} with $\lLie{g}$ replacing $\Lie{g}$.
		\end{enumerate}
\end{proof}

\begin{remark}[The Kostant conjecture about the principal basis]
\label{rem_conj}
Let $\text{Cl}(\Lie{g},K)$ be the Clifford algebra of the Lie algebra $\Lie{g}$ and $K$ its Killing form. Consider the Harish-Chandra projection $\Phi:\text{Cl}(\Lie{g})^{\Lie{g}} \rightarrow \text{Cl}(\Lie{h})$ defined by the decomposition $\text{Cl}(\Lie{g})^{\Lie{g}} = \text{Cl}(\Lie{h}) \oplus \text{Cl}(\Lie{g}) \mathfrak{n}_+ \cap \text{Cl}(\Lie{g})^{\Lie{g}}$, where $\mathfrak{n}_+$ are the positives nilpotent elements. In particular Kostant showed that the image of the Harish-Chandra projection of the primitive generators of $\text{Cl}(\Lie{g})^{\Lie{g}} \cong \bigwedge \Lie{g}^\Lie{g}$ are contained in $\Lie{h} \subset \text{Cl}(\Lie{h})$. Kostant conjectured that the Harish-Chandra projection of the primitive generators of $\Lie{g}$ gives the principal basis of the Cartan Subalgebra of the Langlands dual $\Lie{g}^\vee$ (here using $K$, $\Lie{h}$ and $\lLie{h}$ are identified). This conjecture was partially proven by Bazlov in his thesis \cite{Ba}.
\end{remark}

%---------------------------------------------------------------------------------------------

\section{The $so(2l,\C)$ case ($l$ even)}
\label{sec_d2l}

Throughout this section $\Lie{g}=so(2l,\C)$ with $l$ even, $\Lie{h} \subset \Lie{g}$ is a Cartan subalgebra, $\Weyl$ the Weyl group and $B$ is an invariant bilinear form. The root system of $\Lie{g}$ is $D_l$ and the exponents are given by \cite{Bo}, Table IV

$$\{ 1,3,5,\cdots,l-1,l-1,\cdots,2l-5,2l-3\},$$
Then in the decomposition of $\Lie{g}$ into irreducible $\Lie{s}_0$-modules,

$$\Lie{g} = V_{1} \oplus \cdots \oplus V_{l/2} \oplus V_{l/2+1} \oplus \cdots\oplus V_{l},$$ 
the modules $V_{l/2}$ and $V_{l/2+1}$ have the same dimension $2l-1$.\\

This decomposition does not determine $V_{l/2}$ and $V_{l/2+1}$ uniquely. Indeed $V_{l/2}$ and $V_{l/2+1}$ are generated by the highest weight vectors $v_1$ and $v_2$, i.e. $V_{l/2}=U(\Lie{s}_o)v_1$ and $V_{l/2+1}=U(\Lie{s}_o)v_2$. But $v_1 + v_2$ and $v_1 - v_2$ also generate two irreducible $\Lie{s}_o$-modules whose the direct sum is equal to $V_{l/2} \oplus V_{l/2+1}$, i.e. $U(\Lie{s}_o)(v_1+v_2) \oplus U(\Lie{s}_o)(v_1 - v_2) = V_{l/2} \oplus V_{l/2+1}$.\\ 

The following decomposition is uniquely determined by the action of the principal $\sltwo$-triple,

\begin{equation}
\label{equ_dec_d2l}
	\Lie{g} = V_{1} \oplus \cdots \oplus V_{l/2-1} \oplus W \oplus V_{l/2+2} \oplus \cdots\oplus V_{l},
\end{equation}
where $W=V_{l/2} \oplus V_{l/2+1}$. Note that this decomposition is orthogonal relatively to $B$.\\

Hence the principal basis is uniquely determined modulo a choice of two linearly independent vectors in $W \cap \Lie{h}$, i.e.

$$\text{principal basis of }\Lie{h} = \{h_1, \cdots , h_l , h_{l+1}, \cdots ,h_r \},$$
with $h_i$ such that $<h_i>_\C = V_i \cap \Lie{h}$ for all $i \neq l/2,l/2+1$ and $<h_{l/2},h_{l/2+1}>_\C = W \cap \Lie{h}$. Moreover if we chose $h_{l/2}$ and $h_{l/2+1}$ to be orthogonal, then the principal basis is uniquely determined modulo a rotation in $W$.\\

Now we will give the modifications in the proof of Theorem \ref{thm_main_1}.

\begin{enumerate}
	\item[(a) (b)] The first and second step remain true.
	\item[(c)] The third step is modified as follows: for all $i \neq l/2,l/2+1$ we have
	$$D\tilde{p}_i(\rho^\vee) \in \left( V_{1} \oplus \cdots \oplus V_{i} \right) \cap \Lie{h},$$
	and for $i=l/2$ or $i=l/2+1$ we have
	$$D\tilde{p}_i(\rho^\vee) \in \left( V_{1} \oplus \cdots \oplus W \right) \cap \Lie{h}.$$
	\item[(d)] The fourth step reaming true with a suitable choice of basis vectors in $W \cap \Lie{h}$.
\end{enumerate}

Let $\{\alpha_1, \cdots , \alpha_l\}$ be the $l$ simple roots of $\Lie{g}$. We define the automorphism $\sigma$ by

$$\sigma(\alpha_i)=\alpha_i \quad (i=1,\cdots l-2) \qquad \sigma(\alpha_{l-1})=\sigma(\alpha_l) \qquad \sigma(\alpha_{l})=\sigma(\alpha_{l-1}).$$

This automorphism induces an automorphism on the Lie algebra, see Section 7.9 of \cite{Ka}. As this automorphism is of order two it is diagonalizable and its eigenvalues are $+1$ and $-1$. Let $\Lie{g} = \Lie{g}_0 \oplus \Lie{g}_1$ be the decomposition in the eigenspace $\Lie{g}_0$, resp. $\Lie{g}_1$, of eigenvalues $+1$, resp. $-1$. In \cite{Ka}, Chapter 8, Kac showed that the roots system of $\Lie{g}_0$ is $B_{l-1}$. This implies that the dimension of $\Lie{g}_1$ is $2l-1$. Indeed $\dim(D_{l})=l(2l-1)$ and $\dim(B_{l-1})=l(2l-1)$.\\

We denote by $B$ the canonical form of $\Lie{g}$, i.e. $B(\alpha_i,\alpha_i)=2$.\\

\begin{theorem}
\label{thm_auto}
	Consider the decomposition (\ref{equ_dec_d2l}) of $\Lie{g}$, we have
	
	\begin{enumerate}
		\item	All the irreducible $\Lie{s}_0$-modules $V_i$ are contained in the eigenspace $\Lie{g}_0$.
		\item The $\Lie{s}_0$-module $W$ is invariant by $\sigma$, i.e. $\sigma(W)=W$.
		\item The automorphism $\sigma$ provides the decomposition of $W$ into $W = W_0 \oplus W_1$, where $W_i \subset \Lie{g}_i$. Moreover, this is a decomposition into two irreducible $\Lie{s}_0$-submodules and they are orthogonal.
	\end{enumerate}
	
	Hence if we require that the elements of the principal basis should be eigenvectors of $\sigma$, then the principal basis is uniquely determined (up to constants).
\end{theorem}

The proof of this Theorem follows from the next three propositions.\\

\begin{proposition}
\label{prop_inv1}
	Let $\Lie{s}_0 = <e_0,f_0,h_0>$ be the principal $\sltwo$-triple. Then its adjoint action commute with the automorphism $\sigma$, i.e.
	
	$$\sigma \circ ad_{e_{0}} = ad_{e_{0}} \circ \sigma \qquad \sigma \circ ad_{f_{0}} = ad_{f_{0}} \circ \sigma \qquad \sigma \circ ad_{h_{0}} = ad_{h_{0}} \circ \sigma.$$
\end{proposition}
\begin{proof}
	The principal TDS $\Lie{s}_0$ is given by  
	$$h_0 = 2(l-1)\alpha_1^\vee + 2(2l-3)\alpha_2^\vee + \cdots + (l-2)(l+1)\alpha_{l-2}^\vee + l(l-1)/2(\alpha_{l-1}^\vee + \alpha_{l}^\vee),$$
	$$ e_0 = e_1 + \cdots e_l,$$ 
	$$f_0= 2(l-1)f_1 + \cdots + (l-2)(l+1)f_{l-2} + l(l-1)/2(f_{l-1} + f_{l}),$$
	where $e_i$ is the positive nilpotent element corresponding to the simple root $\alpha_i$, and the $f_i$ is the corresponding negative nilpotent elements such that $B(e_i,f_i)=1$. 
	Clearly they are eigenvectors of eigenvalues $+1$. The result follows from the equation 
	$$\sigma (ad_{e_0} x) = ad_{\sigma(e_{0})} \sigma(x) =  ad_{e_0} \sigma(x).$$
	The proof for $h_0$ and $f_0$ is similar.
\end{proof}
	
\begin{proposition}
\label{prop_inv2}
	Let $V \in \Lie{g}$ be a irreducible $\Lie{s}_0$-module. If there exist $x \in V$ such that $\sigma(x)=x$ (or $\sigma(x)=-x$), then $V \subset \Lie{g}_0$ (or $V \subset \Lie{g}_1$). 
\end{proposition}
\begin{proof}
	Let $x \in V$ be a eigenvector of eigenvalue $\lambda$ for $\sigma$ ($\lambda = \pm 1$). Then using Proposition \ref{prop_inv1}, we have that all vectors $y \in V = U(\Lie{s}_0)x$ are eigenvectors of eigenvalue $\lambda$.
\end{proof}

Before the next proposition we introduce the standard choice for the generators of $\ShW$. \\

Let $\{e^*_1, \cdots e^*_l\}$ be the standard orthonormal basis of $\dLie{h}$. The simple roots are given by

$$\alpha_1 = e^*_1-e^*_2 \qquad \cdots \qquad \alpha_{l-1}=e^*_l-e^*_{l-1} \qquad \alpha_{l}=e^*_{l}+e^*_{l-1}.$$
Recall that on the elements of the basis the canonical bilinear invariant form is given by $B(e_i^*,e_j^*)=\delta_{ij}$.
In term of this basis the standard generators are given by (\cite{Me}, Section 3)

$$p_i = \sum_j (e_j^*)^{2i} \quad i=1,\cdots,l-1 \qquad \text{and} \qquad p_e = e^*_1e^*_2 \dots e^*_l.$$

The half sum of positive roots is given by

$$\rho^\vee = (l-1)(e_{1}) + (l-2)(e_{2}) + \cdots + e_{l-1},$$
where $\{e_1,\cdots,e_l \}$ is the dual basis of $\{e^*_1,\cdots,e^*_l\}$.Note that the two polynomials of degree $l$ are $p_e$ and $p_{l/2}$.\\

The differential of these generators computed at $\rho^\vee$ are, 

\begin{eqnarray*}
\label{equ_basis_d2l}
	Dp_1(\rho^\vee)  = & 2 \rho^\vee \\
	Dp_2(\rho^\vee)  = & 4 \cdot 2^3 ( (l-1)^3e_1 + \cdots + e_{l-1}) \\
	\cdots & \\
	Dp_e(\rho^\vee) = & l!\cdot 2^{l-1}e_l \\
	Dp_{l/2}(\rho^\vee) = & \cdots \\
	\cdots & \\ 
	Dp_{l-1}(\rho^\vee) = & (2l-2)2^{2l-3} ((l-1)^{2l-3} e_1 + \cdots + e_{l-1}) 
\end{eqnarray*}

We remark that $<Dp_1(\rho^\vee), \cdots ,Dp_{l-1}(\rho^\vee) >_\C$ is orthogonal to $<Dp_e(\rho^\vee)>_\C$.\\

From these calculations we deduce the next proposition

\begin{proposition}
\label{prop_inv3}
	We have that $Dp_i(\rho^\vee) \in \Lie{g}_{0}$ for all $i=1,\cdots,l-1$ and $Dp_e(\rho^\vee) \in \Lie{g}_{1}$.\\
\end{proposition}

Now we prove the Theorem \ref{thm_auto}.

\begin{proof}[Proof of Theorem \ref{thm_auto}]
	\begin{enumerate}
		\item We orthogonalize the basis (\ref{equ_basis_d2l}) and obtain the principal basis $\{h_1, \cdots ,h_l\}$. Note that all $h_i$ except $h_{l/2}$ are in the space $<e_1, \cdots  ,e_{l-1}>_\C$ and $h_{l/2} \in <e_{l}>_\C$. By Proposition \ref{prop_inv2}, this proves the first assertion.\\
		\item The second assertion follows from Proposition \ref{prop_inv2} and the fact that $W = U(\Lie{s}_{0})<h_{l/2},h_{l/2+1}>_\C$.
		\item For the third assertion we have that $W_0 = U(\Lie{s}_{0})h_{l/2+1}$ and $W_1 = U(\Lie{s}_{0})h_{l/2}$.
	\end{enumerate}
\end{proof}

%---------------------------------------------------------------------------------------------

\bibliographystyle{amsplain}   

\end{document}